\documentclass[11pt,a4paper]{article}
\usepackage{amsmath}
\usepackage{latexsym}
\font\emailfont=cmtt10
\font\koteret=cmbx9 scaled 1800
\font\tenBbb=msbm10 scaled 1100
\font\sevenBbb=msbm7 scaled 1100 \font\fiveBbb=msbm5 scaled 1100
\newfam\Bbbfam
\textfont\Bbbfam=\tenBbb \scriptfont\Bbbfam=\sevenBbb
\scriptscriptfont\Bbbfam=\fiveBbb
\def\Bbb{\fam\Bbbfam\tenBbb}
\def\R{\Bbb R}

\def\Sn{{\Bbb S}^{n-1}}
\def\norma#1{|\!|#1|\!|}
\def\ip#1#2{\langle #1,#2\rangle}
\def\a{\alpha}

\newcommand{\proof}{\emph{Proof}\/.\quad}
\def\thinsquare{\vcenter{\hrule height 0.3pt \hbox{\vrule width0.3pt
height6pt \kern6.5pt \vrule width 0.3pt} \hrule height 0.3pt}}
\def\sofproof{\hfill$\thinsquare$}
\newtheorem{thrm}{Theorem}[section]
\newtheorem{lem}[thrm]{Lemma}
\newtheorem{prop}[thrm]{Proposition}

\author{Yossi Lonke}
\title{\koteret Derivatives of the $L^p$-cosine transform }
\begin{document}
\date{}
\maketitle
\begin{abstract}The $L^p$-cosine transform of an even, continuous function
$f\in C_e(\Sn)$ is defined by:
$$H(x)=\int_{\Sn}|\ip{x}{\xi}|^pf(\xi)\,d\xi,\quad x\in {\R}^n.$$
It is shown that if $p$ is not an even integer then all partial
derivatives of even order of $H(x)$ up to order $p+1$ (including
$p+1$ if $p$ is an odd integer) exist and are continuous
everywhere in ${\R}^n\backslash\{0\}$. As a result of the
corresponding differentiation formula, we show that if
$f$ is a positive bounded function and $p>1$ then $H^{1/p}$
is a support function of a convex body whose boundary has everywhere positive Gauss-Kronecker curvature.

\end{abstract}
\section{Introduction}
Recent research in convex geometry has repeatedly utilized two important
integral transforms of functions defined on the unit sphere $\Sn$ in ${\R}^n$. These
are the \emph{cosine transform}  and the
\emph{spherical Radon transform}, both acting on $C_e^{\infty}(\Sn)$, the space of infinitely
differentiable even functions on $\Sn$, by:
\begin{eqnarray}
Tf(x)&=&\int_{\Sn}|\ip{x}{\xi}|f(\xi)\,d\xi,\qquad\hbox{(cosine transform)}\nonumber\\
Rf(x)&=&\int_{\Sn\cap x^{\perp}}f(\xi)\,d\xi,\qquad\hbox{(spherical Radon transform)}\nonumber
\end{eqnarray} 
where $\ip{}{}$ denotes the scalar product, $d\xi$ the
spherical Lebesgue measure, and $x^{\perp}$ the $n-1$ dimensional subspace orthogonal
to~$x$. It is well known that $T$ and $R$ are both
continuous bijections of $C_e^{\infty}(\Sn)$ onto itself,
(the topology on $C_e^{\infty}(\Sn)$ taken as uniform convergence of all derivatives).
This fact allows an extension of both transforms, by duality, to bi-continuous bijections
of the dual space $D_e(\Sn)$ of even distributions on $\Sn$. A pleasant consequence of this
extension is that we may assign precise meanings to the symbols $R\rho,R^{-1}\rho, T\rho, T^{-1}\rho$,
for a given even distribution $\rho\in D_e(\Sn)$. For example, one has
$$(T^{-1}\rho)(f)=\rho (T^{-1}f),\qquad\forall\rho\in D_e(\Sn),\quad\forall f\in C_e^{\infty}(\Sn).$$
In particular, one talks about the cosine transform of an $L^1$ function, or the spherical Radon transform
of a measure. These purely analytic manipulations turned out it to have surprisingly far reaching consequences.
For example, the key to the ultimate solution of the Busemann-Petty problem, (which was one of the most intriguing unsolved problems
of convex geometry) was uncovered by Lutwak in \cite{LutwakIntersectionBodies}, where the notion of \emph{intersection body} was
invented. An origin symmetric convex body is called an intersection body if its radial function is realized as
a spherical Radon transform of a positive measure on $\Sn$. Lutwak reduced the Busemann-Petty problem to the
analytic question of whether $R^{-1}\rho$ is a positive measure whenever $\rho$ is a radial function of a centrally
symmetric convex body. The answer is yes, if and only if the dimension is at most~$4$. Although in general
it was known that for sufficiently large $n$ the Busemann-Petty problem has a negative answer in ${\R}^n$
(see \cite{KBall}), the curious dependence
on the dimension and the precise role of convexity
were not understood until they were revealed by means of sophisticated analysis in \cite{GKS}.

The relevance of the cosine transform to convex geometry becomes clear through the concept of \emph{zonoids},
also called \emph{projection bodies}.
These are bodies that can be approximated to any degree of accuracy, in the Hausdorf metric sense, by finite vector sums of intervals, called
\emph{zonotopes}. Every zonotope has a center of symmetry (namely, the sum of the centers of the intervals).
Up to translation, every zonotope $Z$ has therefore the form $Z=\sum_1^m\lambda_i[-u_i,u_i]$, for some
positive numbers $\lambda_i$ and $u_i\in\Sn$. Here $[-u_i,u_i]$ denotes the convex hull of $\{-u_i,u_i\}$. The support function of $Z$ is then $h_z(x)=\sum_1^m\lambda_i|\ip{u_i}{x}|$.
Let $\delta_u$ denote the unit-mass measure concentrated at $u\in\Sn$. Put ${\mu=\sum_1^m\lambda_i\frac{\delta_{u_i}+\delta_{-u_i}}{2}}$.
Then
\begin{equation}\label{zonotope}
h_Z(x)=\int_{\Sn}|\ip{x}{u}|\,d\mu=T\mu(x).
\end{equation}
In other words, the support function of a zonotope is a cosine transform of a positive, discrete measure. A standard
approximation argument yields a fundamental theorem: \emph{A centrally symmetric convex body is a zonoid if and only if
its support function is a cosine transform of a positive measure}.

The measure $\mu$ in (\ref{zonotope}) is called the \emph{generating measure} of~$Z$. Generalizing this concept,
Weil \cite{Weildis} proved that to every centrally symmetric convex body $K\subset{\R}^n$ corresponds a unique \emph{generating distribution},
that is, a continuous linear functional $\rho_K$ on the space $C_e^{\infty}(\Sn)$, whose domain can be extended as to 
include the functions $|\ip{u}{\cdot}|$ with $u\in\Sn$, such that $\rho_K(|\ip{u}{\cdot}|)=h_K(u)$ for every $u\in\Sn$.
Recall that positive distributions are in fact positive measures. Thus in the context of zonoids
Weil's result is particularly useful --- it provides \emph{a-priori} a functional, namely $T^{-1}h_K$, whose positivity
is to be checked. Interestingly, the cosine and spherical Radon transforms are related by:
\begin{equation}\label{inverse}
T^{-1}=c_n(\Delta_n+n-1)R^{-1},
\end{equation}
where $\Delta_n$ is the spherical Laplace operator on $\Sn$, and $c_n>0$ (see \cite{GoodWeil}).
The inversion formula (\ref{inverse}) proved a useful analytic tool in
constructing examples of non-smooth zonoids whose polars are zonoids \cite{Lonkezon}, and
of convex bodies whose generating distributions have
large degree \cite{LonkeGenDist}. 
 
Often one thinks of $h_Z(x)$ in (\ref{zonotope}) as representing the
norm of some space, which in this case is isometric to a subspace
of $L^1(\Sn,\mu)$. A natural generalization is then to look at
functions of the form
\begin{equation}\label{q-reprsnt}
H^p(x)=\int_{\Sn}|\ip{x}{\xi}|^p\,d\mu,\quad (p\geq 1)
\end{equation}
If $\mu$ is positive, $H$ is continuous, convex and
$1$-homogeneous, hence a support function of some convex body,
and also the norm of some normed space, which is evidently
isometric to a subspace of $L^p(\Sn,\mu)$. The r.h.s of~(\ref{q-reprsnt})
 is called the \emph{$L^p$-cosine transform} of the measure $\mu$, and
is denoted by $T_p\mu$.
If $p$ is not an even integer, the measure $\mu$ in
(\ref{q-reprsnt}) is uniquely determined by the norm on the left
hand side. For $p=1$, this was first proved by Alexandrov
\cite{Alexandrov} and rediscovered several times since. In
\cite{Neyman}, Neyman proved that if $p$ is not an even integer,
the linear span of the functions $|\ip{x}{\cdot}|^p$, defined on
$\Sn$ and indexed by $x\in{\R}^n$, is dense in the space
$C_e(\Sn)$ of continuous even functions on $\Sn$. In particular,
$\mu$ in (\ref{q-reprsnt}) is uniquely determined. If $p$ is an even integer, the
functions $|\ip{x}{\cdot}|^p$ span precisely the subspace of
homogeneous (even) polynomials of degree $p$ (see \cite{Neyman}), so
that there is no longer uniqueness in the representation
(\ref{q-reprsnt}). The inversion problem for the $L^p$-cosine transform
of $L^1$ functions
has been treated in \cite{KoldHouston} in several important special cases. The general case
of inversion has apparently been neglected.

In a recent paper \cite{frenchguy}, the cosine transform of
a continuous function was shown to be a $C^2$ function. In the first section below, this
result is generalized in two ways. First, it is proved that for a nonnegative integer~$k$,
the~${2k+1}$-cosine transform of a continuous function
is of class $C^{2k+2}$. The proof below invokes Fourier transform techniques developed by
Koldobksy in a series of papers (\cite{Kold4,KoldHouston,Kold3,Kold1,Kold2}).
Then, we deal with the $L^p$-cosine transform where $p>1$ is not an integer,
and show that if $f$ is a bounded function, then $T_pf$ has continuous
derivatives of the largest even order smaller than $p+1$. For second order derivatives,
this was done in a more general setting in \cite{KoldIndiana}, using other methods. The first section
is concluded with an additional result, asserting that the cosine transform carries $L^1(\Sn)$ into
$C^1(\Sn)$. Our main application
is expounded in section~2, where we show that if $H^p=T_pf$ with $f$ positive and bounded, then $H$ is a support function of a centrally symmetric $C_+^2$ convex body.
That is, the boundary of the body has everywhere positive Gauss-Kronecker curvature. This should be compared
to Theorem~2 of \cite{frenchguy}, which asserts that zonoids (i.e, the $p=1$ case) whose generating measures
are continuous functions may fail to have positive Gauss-Kronecker curvature at some boundary point only
because all the principal radii of curvature evaluated at the corresponding outward unit normal
are zero (whereas in general the curvature may not exist due to just one vanishing principal radius of curvature).

\section{Differentiation of the $L^p$-cosine transform}
Let $\alpha=(\alpha_1,\dots,\alpha_n)$ denote
a multi-index ( $\alpha_k$ are nonnegative integers ).
${|\alpha|=\alpha_1+\cdots+\alpha_n}$. Given $\alpha$,
$D^{\alpha}$ denotes the differential operator
$$D^{\alpha}=(\frac{\partial}{\partial x_1})^{\alpha_1}
(\frac{\partial}{\partial x_2})^{\alpha_2}\cdots(\frac{\partial}{\partial x_n})^{\alpha_n}$$
In what follows, we denote
$C_t=\frac{2^{t+1}\sqrt{\pi}\Gamma((t+1)/2)}{\Gamma(-t/2)}$. Our first result
is a generalization of Th.~1 in \cite{frenchguy}.
\begin{thrm} Let $n\geq 2$ and suppose that
$$
H(x)=\int_{S^{n-1}}|\ip{x}{\xi}|^{2k+1}\,f(\xi)\,d\xi
$$
where $k$ is a nonnegative integer, and $f\in C_e(\Sn)$. Then $H\in C^{2k+2}({\R}^n\backslash\{0\})$
and for every multi-index $\a$ with $|\a|=2k+2$, one has for each $x\in{\R}^n\backslash\{0\}$
\begin{equation}\label{q-odd-int}
D^{\a}H(x)=C_{2k+1}\frac{(-1)^{k+1}}{\norma{x}}\int_{\Sn\cap
x^{\perp}}\xi_1^{\alpha_1}\cdots\xi_n^{\alpha_n}f(\xi)\,d\xi
\end{equation}
\end{thrm}
\vskip 6pt
In case the differentiation-order $|\a|$ strictly smaller than $p+1$, (and even) the assumptions on $f$ can be somewhat relaxed, and the corresponding
differentiation formula is different. For these reasons the result is formulated separately.

\begin{thrm}
Let $n\geq 2$ and suppose that
$$
H(x)=\int_{S^{n-1}}|\ip{x}{\xi}|^p\,f(\xi)\,d\xi
$$
where $p>1, p\neq 2k$ and $f\in L^{\infty}(\Sn)$. Let
$\a$ be a multi-index such that $|\a|$ is even and $|\a|<p+1$.
Then $H\in C^{|\alpha|}({\R}^n\backslash\{0\})$ and
\begin{equation}\label{q_non_integer}
D^{\alpha}H(x)=i^{|\alpha|}\frac{C_p}{C_{p-|\a|}}\int_{\Sn}|\ip{x}{\xi}|^{p-
|\alpha|}\xi_1^{\alpha_1}\cdots\xi_n^{\alpha_n}f(\xi)\,d\xi.
\end{equation}
\end{thrm}
\vskip 12pt
For the proofs, we use distribution theory and Fourier transforms. As usual, let
$S({\R}^n)$ denote the space of rapidly decreasing infinitely differentiable functions (test functions)
in ${\R}^n$, and $S^{'}({\R}^n)$ is the space of distributions over $S({\R}^n)$. The Fourier
transform of a distribution $f\in S^{'}({\R}^n)$ is defined by $({\hat{f}},{\hat{\phi}})=(2\pi)^n({f},{\phi})$,
for every test function $\phi$. 
\vskip 12pt
\emph{Proof of Theorem~2.1}\quad
For every test function $\phi(x)$ supported in ${\R}^n\backslash\{0\}$ consider another test function
$\psi(x)=x_1^{\a_1}\cdots x_n^{\a_n}\phi(x)$. Since $|\a|$ is
even, so is $\psi$. From lemma 2.2 of \cite{KoldHouston} we have
\begin{equation}\label{firstcase}
(\hat{H},\psi)=(2\pi)^{n-1}C_{2k+1}\int_{\Sn}f(\xi)\,d\xi\int_{\R}t^{-2k-2}\psi(t\xi)\,dt.
\end{equation}

Therefore,
$$
(\prod_{k=1}^nx_k^{\a_k}\hat{H},\phi)=(2\pi)^{n-1}C_{2k+1}\int_{\Sn}\prod_{k=1}^n
\xi_k^{\a_k}f(\xi)\,d\xi\int_{\R}\phi(t\xi)\,dt.
$$
By the well-known connection between the Fourier transform and the Radon transform (see \cite{GelfandShilov}),
the function $t\to(2\pi)^n\phi(-t\xi)$ is the Fourier transform of the function
$z\to\int_{\ip{x}{\xi}=z}\hat{\phi}(x)\,dx$. Therefore,
$\int_{\R}\phi(t\xi)\,dt=(2\pi)^{-n+1}\int_{\xi^{\perp}}\hat{\phi}(x)\,dx$, so
we have
$$
(\prod_{k=1}^nx_k^{\a_k}\hat{H},\phi)=C_{2k+1}\int_{\Sn}\prod_{k=1}^n
\xi_k^{\a_k}f(\xi)\,d\xi\int_{\xi^{\perp}}\hat{\phi}(x)\,dx
$$
Put $g(\xi)=\prod_{k=1}^n\xi_k^{\a_k}f(\xi)$, and let $R$ denote
the spherical Radon transform. Since for $n\geq 2$ the
function $\norma{x}_2^{-1}$ is locally integrable, we have
$$\aligned
\int_{{\R}^n}\norma{x}_2^{-1}\hat\phi(x)Rg(x/\norma{x}_2)\,dx &=
\int_0^{\infty}r^{n-2}\left(\int_{\Sn}\hat{\phi}(r\xi)Rg(\xi)\,d\xi\right)\,
dr\\
&=\int_0^{\infty}r^{n-2}\int_{\Sn}g(\xi)\,d\xi\int_{\xi^{\perp}\cap\Sn}\hat{
\phi}(ru)\,du\,dr\\
&=\int_{\Sn}g(\xi)\,d\xi\int_{\xi^{\perp}}\hat{\phi}(x)\,dx
\endaligned
$$
Self-duality of the spherical Radon transform was used here.
Consequently,
\begin{equation}\label{first}
(\prod_{k=1}^nx_k^{\a_k}\hat{H},\phi)=C_{2k+1}\int_{{\R}^n}\norma{x}_2^{-1}\hat\phi(x)Rg(x/\norma{x}_2)\,dx
\end{equation}

On the other hand, the well known connection between
differentiation and Fourier transforms gives:
\begin{equation}\label{second}
(\prod_{k=1}^nx_k^{\a_k}\hat{H},\phi)=i^{-|\a|}(D^{\alpha}H,\hat{\phi})
\end{equation}
Recall that $\hat{\hat{\phi}}=(2\pi)^n\phi(-x)$. Therefore, for every distribution $f$ and
an even test function $\phi$, one has $(f,\hat{\phi})=(\hat{f},\phi)$. 
Since $\phi(x)$ is an arbitrary even test function (with $0\notin$
supp$\,\phi$) (\ref{first}), (\ref{second})Ê together imply that the
Fourier transforms of the distributions
\begin{equation}\label{third}
D^{\alpha}H(x)\quad\text{and}\quad
C_{2k+1}\frac{(-1)^{k+1}}{\norma{x}_2}Rg(\frac{x}{\norma{x}_2})
\end{equation}
are equal distributions in ${\R}^n\backslash\{0\}$. Therefore, the distributions in
(\ref{third}) can differ by a polynomial only (\cite{GelfandShilov}, p. 119). Since both distributions
are even and homogeneous of degree $-1$, the polynomial
must be identically zero. Hence the distributions in (\ref{third}) are equal. To show that $H$ is a $C^{|\a|}$
function we must show that $D^{\a}H$ exists also in the classical sense and is continuous.
As is well known in the theory of distributions, classical and distributional derivatives coincide if
the distributional derivative in question happens to be a continuous function. (\cite{LiebLoss}, p. 136).
Since $f$ is continuous, so is the spherical Radon transform of $\prod_{k=1}^n\xi_k^{\a_k}f(\xi)$.
Therefore $D^{\a}H(x)$ is a continuous function, and we have (\ref{q-odd-int}).\sofproof\vskip 6pt
The proof of Theorem~2.2 uses the same technique.
Instead of (\ref{firstcase}) we now have:
\begin{equation}\label{qcase}
(\hat{H},\psi)=(2\pi)^{n-1}C_p\int_{\Sn}f(\xi)\,d\xi\int_{\R}|t|^{-1-p}\psi(t\xi)\,dt
\end{equation}

Therefore,
$$
(\prod_{k=1}^nx_k^{\a_k}\hat{H},\phi)=(2\pi)^{n-1}C_p\int_{\Sn}\prod_{k=1}^n
\xi_k^{\a_k}f(\xi)\,d\xi\int_{\R}|t|^{|\a|-p-1}\phi(t\xi)\,dt.
$$
Since $p-|\a|>-1$, and $p-|\a|$ is not an even integer, we can
apply Lemma 2.1 of \cite{KoldHouston}:
$$\int_{\R}|t|^{|\a|-p-1}\phi(t\xi)\,dt=\frac{1}{(2\pi)^{n-1}C_{p-|\a|}}
\int_{{\R}^n}|\ip{x}{\xi}|^{p-|\a|}\hat{\phi}(x)\,dx$$

Consequently,
\begin{equation}\label{qfirst}
(\prod_{k=1}^nx_k^{\a_k}\hat{H},\phi)=\frac{C_p}{C_{p-|\a|}}\int_{{\R}^n}
\left[\int_{\Sn}|\ip{x}{\xi}|^{p-|\a|}\prod_{k=1}^n\xi_k^{\a_k}f(\xi)\,d\xi\right]\hat{\phi}(x)\,dx
\end{equation}
The connection between differentiation and the Fourier transform yields in this case: 
\begin{equation}\label{qsecond}
(\prod_{k=1}^nx_k^{\a_k}\hat{H},\phi)=i^{-|\a|}(D^{\alpha}H,\hat{\phi})
\end{equation}
Together, (\ref{qfirst}) and (\ref{qsecond}) imply that the
Fourier transforms of the distributions
\begin{equation}\label{qthird}
D^{\alpha}H(x)\quad\text{and}\quad
i^{|\a|}\frac{C_p}{C_{p-|\a|}}\int_{\Sn}|\ip{x}{\xi}|^{p-|\a|
}\prod_{k=1}^n\xi_k^{\a_k}f(\xi)\,d\xi
\end{equation}
are equal distributions in ${\R}^n\backslash\{0\}$. As before, 
the distributions in (\ref{qthird}) are equal.  It
remains to check that the right hand side of (\ref{qthird}) is a continuous function. 
It is obviously continuous in
$x$ if $p>|\a|$. To see that it is also continuous in the case
$|\a|-1<p<|\a|$, pick a sequence $x_m\neq 0$ such that
$\lim_{m\to\infty}x_m=x_0\neq 0$. For sufficiently large $m$, we
have $|\ip{x_m}{\xi}|\geq |\ip{x_0}{\xi}|/2$ for each $\xi\in\Sn$.
Therefore, the integrand in the right hand side of (\ref{qthird}) is almost
everywhere bounded above by the function $\xi\to
(|\ip{x_0}{\xi}|/2)^{p-|\a|}\norma{f}_{\infty}$, which is in
$L^1(\Sn)$ since $p-|\a|>-1$. The desired continuity now follows from
Lebesgue's bounded convergence theorem. \sofproof\vskip 6pt
We conclude this section with a supplementary result, related to the $k=0$ case of Th.~2.1 above.
\begin{prop} 
The cosine transform of an $L^1$ function is continuously differentiable in ${\R}^n\backslash\{0\}$.
\end{prop}
\proof Let $f\in L^1$. By linearity of $T$ and by writing $f=f_+-f_-$, where $f_+,f_-$ are the positive and negative parts of $f$
respectively, we may assume $f\geq 0$. In that case, $Tf$ is a support function of a zonoid~$Z$. Put $h_Z=Tf$. A convex body is
strictly convex (i.e., contains no line-segments in its boundary) if and only if its support function differentiable
in ${\R}^n\backslash\{0\}$ (\cite{Schneider}, 1.7.3, p.~40). If a zonoid $Z$ is not strictly convex, its boundary
has some lower dimensional face, which must be a translate of a zonoid of lower dimension that is a summand
of~$Z$. (\cite{Bolk}, Th.~3.2). This means that $Z$ can be decomposed as $Z=Z_1+Z_2$ where at least one of the summands has lower dimension.
It follows that the generating measure of $Z$ is a sum $\mu_1+\mu_2$ of the generating measures of 
$Z_1,Z_2$, and at least one of these measures is supported on a proper subspace. In particular, $\mu_1+\mu_2$
is not absolutely continuous; but the generating measure of $Z$ is. Thus if $h_Z=Tf$ with $f\in L^1$ and $f>0$,
then $Z$ is a strictly convex zonoid, so $h_Z$ is differentiable in ${\R}^n\backslash\{0\}$. The proof is completed
by noting that support functions differentiable in ${\R}^n\backslash\{0\}$ are already continuously differentiable
there.\sofproof
\section{Application to curvature and convexity}
The main result in this section is the following
\begin{thrm}
Suppose $n\geq 2$ and
$$
H^p(x)=\int_{S^{n-1}}|\ip{x}{\xi}|^p\,f(\xi)\,d\xi
$$
where $p>1, p\neq 2k$ and $f$ is a positive element of $L^{\infty}(\Sn)$. Then $H(x)$ is 
a support function of a centrally symmetric convex body that has everywhere positive
Gauss-Kronecker curvature.
\end{thrm}
The proof is largely based upon the next lemma.\begin{lem} Assume $H^p=T_pf$, where $p$ and $f$ are as in Theorem~3.1. For every unit vector
$u\in\Sn$, and every $v\neq u$, the second directional derivative of $H$ in the direction of $u$, evaluated
at $v$, is positive.
\end{lem}
\proof Differentiating, one finds:
\begin{equation}\label{derivatives}
\frac{\partial^2H}{\partial x_i\partial x_j}=\frac{1}{H^{p-1}}\frac{1}{p}
\left[\frac{\partial^2H^p}{\partial x_i\partial x_j}-\frac{p-1}{p}\frac{1}{H^p}\frac{\partial H^p}{\partial x_i}\frac{\partial H^p}{\partial x_j}\right]
\end{equation}
If $H^p=T_pf$, then by Theorem 2.2
$$\frac{\partial^2H^p}{\partial x_1^2}=p(p-1)\int_{\Sn}|\ip{u}{\xi}|^{p-2}\xi_1^2f(\xi)\,d\xi$$
Moreover, differentiation under the integral sign can easily be justified and
$$\frac{\partial H^p}{\partial x_1}=p\int_{\Sn}|\ip{u}{\xi}|^{p-1}{\rm sgn}\ip{u}{\xi}\xi_1f(\xi)\,d\xi$$
Next, applying the triangle inequality and the Cauchy-Schwartz inequality:
$$
\aligned
&\left|\int_{\Sn}|\ip{u}{\xi}|^{p-1}{\rm sgn}\ip{u}{\xi}\xi_1f(\xi)\,d\xi\right|\leq\int_{\Sn}|\ip{u}{\xi}|^{p-1}|\xi_1|f(\xi)\,d\xi\\
&\leq\left(\int_{\Sn}|\ip{u}{\xi}|^{p-2}\xi_1^2f(\xi)\,d\xi\right)^{1/2}\left(\int_{\Sn}|\ip{u}{\xi}|^{p}f(\xi)\right)^{1/2}\\
&=\left(\frac{1}{p(p-1)}\frac{\partial^2H}{\partial x_1^2}\right)^{1/2}H^{p/2}\\
\endaligned
$$
Therefore,
$$
\left(\frac{\partial H^p}{\partial x_1}\right)^2\leq\frac{p}{p-1}\frac{\partial^2H^p}{\partial x_1^2}H^p,
$$
which implies $\frac{\partial^2H}{\partial x_1^2}\geq 0$. In case of equality, we have equality in
the triangle inequality, and in the Cauchy-Schwartz inequality, applied to the functions $|\ip{u}{\cdot}|^{\frac{p-2}{2}}|\xi_1|$ and
$|\ip{u}{\cdot}|^{\frac{p}{2}}$. Therefore, for every $\xi\in{\rm supp}\,f$, we have for some
real constants $s,t$ not both zero:
\begin{eqnarray}
{\rm (i)} \quad{\rm sgn}\ip{u}{\xi}\xi_1&=&|\xi_1|,\nonumber\\
{\rm (ii)} \quad s|\ip{u}{\xi}|^{p-2}\xi^2 & = & t|\ip{u}{\xi}|^{p}\nonumber
\end{eqnarray}
(i) implies $s\xi^2=t\ip{u}{\xi}^2$. We can not have
$s=0$ (resp. $t=0$), for then the support of $f$ would have to be contained in $u^{\perp}$ (resp. $e_1^{\perp}$), which contradicts $\int_{\Sn}f\,d\xi>0$.
Hence both $s,t$ are non zero, and have the same sign, so that with $\lambda=(s/t)^{1/2}$ we have
$\lambda|\ip{u}{\xi}|=|\ip{e_1}{\xi}|$, and we can drop the absolute values, because $\ip{u}{\xi},\ip{e_1}{\xi}$
have the same sign. Consequently,
$$\ip{\xi}{\lambda e_1-u}=0\qquad\forall\xi\in{\rm supp}\,f,$$
which unless $u=e_1$, contradicts the fact that $\int_{\Sn}f\,d\xi>0$.
Therefore, unless $u=e_1$, one has $\frac{\partial^2H}{\partial x_1^2}> 0$.

Now let $u$ be any direction, and let $U$ be an orthonormal matrix such that $Ue_1=u$. Let $D_u$ denote differentiation
in the $u$ direction. A simple calculation yields:
$$D_u(D_uH)(Uv)=\frac{\partial^2 H\circ U}{\partial x_1^2}(v).$$
Since $H^p=T_pf$, one has $(H\circ U)^p=T_p(f\circ U)$. By the first part of the proof, applied to $H\circ U$
and $f\circ U$
in place of $H$ and$ f$, we get: $\frac{\partial^2 H\circ U}{\partial x_1^2}(v)>0$ whenever $v\in\Sn$ and $v\neq e_1$.
Therefore $D_u^2H(v)>0$ whenever $\xi\neq u$, as was asserted.
\sofproof

\emph{Proof of Theorem 3.1}\quad By Theorem 2.2, $H^p$, and therefore $H$, are $C^2$ functions in ${\R}^n\backslash\{0\}$ Since $f$ is positive,
$H$ is a support function of some (strictly) convex body, say, $K$. To show that $K$ is of class $C_+^2$, it suffices to
show that $K$ has everywhere positive principal radii of curvature. (\cite{Schneider}, p.~111). Let $T_u$ denote
the tangent space to $\Sn$ at $u$. The principal radii of curvature
are eigenvalues of the \emph{reverse Weingarten map} $\overline{W}_u:T_u\to T_u$, where $\overline{W}_u$ is $d(\nabla H)_u$.
( Note that since the gradient $\nabla H(u)$ is the unique point on the boundary of $K$ at which $u$ is an outer normal
vector, its gradient $d(\nabla H)_u$ maps the tangent space $T_u$ into itself). By \cite{Schneider},  p.~108, Lemma~2.5.1,
$$d^2H_u(v,w)=\ip{\overline{W}_uv}{w},\qquad\forall v,w\in T_u.$$
Therefore, if $\lambda$ is an eigenvalue of $\overline{W}_u$ with an eigenvector $v$, then $\lambda=d^2H_u(v,v)$.
As explained in \cite{Schneider} p.~110, $d^2H_u(v,v)=D_u^2H(v)$, which by Lemma~3.2 is positive. \sofproof\vskip 6pt
{\bf Remark}\quad Theorem~3.1 no longer holds for $p=1$. In fact,
we can have $h=Tf$ with $f\in C_e^{\infty}(\Sn)$ and $f>0$, but nonetheless $h$ is not $C^2_+$. Any zonoid
whose support function is $C^{\infty}$ but not $C^2_+$ will do.
\vskip 6pt
A special case of Theorem~2.1, for $k=0$, was proved (in an elementary way)
recently in \cite{frenchguy}. Clearly, the $L^p$-cosine transform of a positive measure is a convex function,
if $p\geq 1$. However, there are also $L^p$-cosine transforms of signed measures, possibly not positive, that are convex
functions. A~theorem by Lindquist \cite{Lindquist} asserts that the cosine transform $Tf(x)$ defines
a support function if and only if 
\begin{equation}\label{Lindquist}
\int_{\Sn\cap u^{\perp}}\ip{\xi}{x}^2f(\xi)\,d\xi\geq 0
\end{equation}
for all $u\in\Sn$ and all $x\in\Sn\cap u^{\perp}$. As was observed in \cite{frenchguy},
the expression in (\ref{Lindquist}) is precisely $d^2H_u(x,x)$, where $H=Tf(x)$. Thus, \emph{a-posteriori} Lindquist's
criterion reduces to the classical assertion that a positively $1$-homogeneous function (i.e, $Tf$ ) is a support
function if and only if its second differential is positive semidefinite at every point.
In this case, homogeneity of $H$ permits consideration of $d^2H_u(x,x)$
only for $x\perp u$.

Put $k=0$ in (\ref{q-odd-int}). The result is:
\begin{equation}
\frac{\partial^2H}{\partial x_i\partial
x_j}(u)=\frac{2}{\norma{u}_2}\int_{\Sn\cap
x^{\perp}}\xi_i\xi_jf(\xi)\,d\xi\qquad (u\in
{\R}^n\backslash\{0\})
\end{equation}
This in turn implies that for $u\in\Sn$, the Hessian matrix $H^{''}$ evaluated at $u$ is given by:
\begin{equation}
\ip{H^{''}_ux}{y}=2\int_{\Sn\cap
u^{\perp}}\ip{x}{\xi}\ip{y}{\xi}f(\xi)\,d\xi
\end{equation}
Therefore $\ip{H^{''}_ux}{x}$ becomes the integral in (\ref{Lindquist}). All this
was pointed out in \cite{frenchguy}. Applying the same reasoning to (\ref{q_non_integer}), we get for $p>1$
($p$ not an even integer):
\begin{equation}
\frac{\partial^2H}{\partial x_i\partial x_j}(u)=
p(p-1)\int_{\Sn}|\ip{u}{\xi}|^{p-2}\xi_i\xi_jf(\xi)\,d\xi
\end{equation}
Hence we derive the following result -- a $p$-version of Lindquist's criterion, which
is an immediate consequence of the previous equation.
\begin{thrm} Suppose $n\geq 2$ and
$$
H(x)=\int_{S^{n-1}}|\ip{x}{\xi}|^p\,f(\xi)\,d\xi
$$
where $p>1, p\neq 2k$ and $f\in C_e(\Sn)$. Then $H(x)$ is convex
if and only if for all $u\in\Sn, x\in{\R}^n$
\begin{equation}
\int_{\Sn}|\ip{u}{\xi}|^{p-2}\ip{x}{\xi}^2f(\xi)\,d\xi\geq 0
\end{equation}
\end{thrm}
\vskip 12pt
{\bf Acknowledgment} I thank Alexander Koldobsky for enlightening discussions.

\providecommand{\bysame}{\leavevmode\hbox
to3em{\hrulefill}\thinspace}

\vskip 12pt

Email address: {\emailfont yossil1@mac.com}

\end{document}